\documentclass[12pt]{article}
\usepackage{amsmath}
\usepackage{amscd}
\usepackage{mathrsfs}
\usepackage{amsfonts}
\usepackage{amsmath,amssymb}
\baselineskip 5.5pt \lineskip 5.5pt \numberwithin{equation}{section}

\def \Z{\hbox{$Z\hskip -5.2pt Z$}}

\def \C{\hbox{$C\hskip -5pt \vrule height 6pt depth 0pt \hskip 6pt$}}

\def\qed{\ \ \ifhmode\unskip\nobreak\fi\ifmmode\ifinner
         \else\hskip5pt\fi\fi
 \hbox{\hskip5pt\vrule width4pt height6pt depth1.5pt\hskip 1 pt}}
\def\a{\alpha}

\def\d{\delta}

\def\cl{\centerline}

\def\vs{\vspace*}

\def\C{\mathbb{C}}

\def\Z{\mathbb{Z}}

\newtheorem{theo}{Theorem}[section]
\newtheorem{lemm}[theo]{Lemma}

\newtheorem{prop}[theo]{Proposition}

\begin{document}

\begin{center}
{\Large \textbf{Two-parameter Quantum group coming from two-parameter deformed Virasoro Algebra of Hom-type}
\noindent\footnote{Supported by the National Science Foundation of
China (Nos. 11047030 and 11771122).
}} \vs{6pt}
\end{center}

\cl{Wen Zhou, Yongsheng Cheng}

\cl{ \small School
of Mathematics and Statistics, Henan
University, Kaifeng 475004, China} \vs{6pt}

\vs{6pt}

\noindent{\it \textbf{Abstract.}~}
In this paper, firstly, we use the bosonic oscillators to construct a two-parameter deformed Virasoro algebra, which is a non-multiplicative Hom-Lie algebra.
Secondly, a non-trivial Hopf structure related to the two-parameter deformed Virasoro algebra is presented, that is, we construct a new two-parameter quantum group.

{\small
\parskip .005 truein
\baselineskip 10pt \lineskip 10pt
\cl{\bf\S1. \ Introduction}
\setcounter{section}{1}\setcounter{theo}{0}\setcounter{equation}{0}
The $q$-deformed Virasoro algebras are given by many authors (cf. \cite{[ACKP],
[ELMS], [HLS], [Hu2], [K], [OS], [Yau]}) respectively, which can be viewed as a typical examples of the physical applications of quantum group.
Quantum groups are a kind of non-commutative and cocommutative Hopf algebras, which were introduced
by Drinfeld and Jimbo as a $q$-deformation of the universal enveloping algebra of a Lie algebra \cite{[D], [DG], [FR], [MT]}.
Two-parameter quantum deformation is a generalization of the one-parameter quantum deformation. Two-parameter quantum enveloping algebras are known to have
a generalized root space structure and the Drinfeld realizations of the two-parameter quantum enveloping algebras were studied in \cite{[HRZ]}.

As a generalization of Lie algebras, Hom-Lie algebras were
introduced by Hartwig, Larsson and Silvestrov in \cite{[HLS]} as part of a
study of deformations of the Witt and the Virasoro algebras.
The motivations to study Hom-Lie structures are related to physics and
to deformations of Lie algebras, in particular Lie algebras of vector field \cite{[CS1], [CS2], [ELMS], [Yau]}.
A {\it Hom-Lie algebra} is a triple $(L, [\cdot,\cdot], \a)$, in
which $L$ is a vector space, $\a$ is an endomorphism of $L$, and the
skew-symmetric bracket satisfies the following conditions\vs{-6pt}
\begin{equation}\label{Lie-s-s}  [x,y]=-[y,x]
\ \mbox{ \ ({\it skew symmetry}),}
\end{equation}
\begin{equation}\label{Lie-j-i}
[\a(x),[y,z]]+[\a(y),[z,x]]+[\a(z),[x,y]]=0, \ \forall x,y,z\in V \
\mbox{ \ ({\it generalized Jacobi identity}). }
\end{equation}
Obviously, Lie algebras are special cases of Hom-Lie algebras in
which $\a$ is the identity map.

In \cite{[ELMS]}, O. Elchinger et. al introduced the two parameters deformed Virasoro algebra $V_{p,q}$, which is a Hom-Lie algebra.
$V_{p,q}=(\hat{L},\hat{\alpha})$ has basis $\{L_n,C|n\in\Z\}$ and bracket relations:
\begin{align*}
[L_n,L_m]:&=(\frac{[n]}{p^{n}}-\frac{[m]}{p^{m}})L_{n+m}
+\d_{m+n,0}\frac{(q/p)^{-n}}{6(1+(q/p)^{n})}
\frac{[n-1]}{p^{n-1}}
\frac{[n]}{p^{n}}\frac{[n+1]}{p^{n+1}}C,\\
[\hat{L},C]:&=0,
\end{align*}
and $\hat{\alpha}:\hat{L}\longrightarrow \hat{L}$ is the endomorphism of $\hat{L}$ defined by
$\hat{\alpha}(L_n)=((1+(q/p)^{n}))L_n,\ \hat{\alpha}(C)=C$. The main tools are $(\sigma, \tau)$-derivations which are generalized
derivations twisting the Leibniz rule by means of a linear map.

In 1998, Hu gave the quantum group structure of the
q-deformed Virasoro algebra in \cite{[Hu1]}. In \cite{[CS1]}, Cheng and Su developed an approach to construct a
q-deformed Heisenberg-Virasoro algebra, which is a Hom-Lie algebra, and the
quantun deformations of Heisenberg-Virasoro algebra which provided a
nontrivial Hopf structure were presented. In \cite{[Yuan]},
Yuan realized the q-deformation $W(2,2)$ by using the bosonic and fermionic oscillators in physics,
the quantum group structure of q-deformation on Lie algebra $W(2,2)$ is further determined.
For the superversion, a two-parameter quantum deformation of Lie superalgebra in the non-standard simple root system with two odd simple roots is constructed in \cite{[HM]}. \par

In the oscillator, the bosonic oscillator $a$ and its hermitian conjugate $a^{+}$ obey the commutation relations:
\begin{align}
&[a,a^{+}]=aa^{+}-a^{+}a=1,\ [1,a^{+}]=[1,a]=0.\label{9.100}
\end{align}
According to \cite{[CS1]}, the Hopf structure on a algebra is as follows.
By a Hopf structure on a algebra $A$, we mean that $A$ is associated with a triple $(\Delta,\epsilon,S)$, where the coproduct $\Delta$: $A\rightarrow A\otimes A$ is an algebra homomorphism, the counit $\epsilon$: $A\rightarrow \mathbb{F}$ is an algenra homomorphism, and the antipode $S$:
$A\rightarrow A$ is an anti-homomorphism such that
\begin{align*}
(1\otimes \Delta)\Delta(x)=(\Delta\otimes 1)\Delta(x)\ \ (coassociativity),\\
m((1\otimes \epsilon)\Delta(x))=x=m((\epsilon \otimes1)\Delta(x))\ \ (counit\ axiom),\\
m((1\otimes S)\Delta(x))=\epsilon(x)=m((S \otimes1)\Delta(x))\ \ (antipode\ axiom),
\end{align*}
for all $x\in A$, where m: $A\otimes A\rightarrow A$ is the multiplication map of A. A Hopf algebra is an algebra equipped with a Hopf structure.

The goal of this paper is twofold. Firstly, we use the bosonic oscillators to construct a two-parameter deformed Virasoro algebra, which is a non-multiplicative Hom-Lie algebra.
Secondly, a non-trivial Hopf structure related to the two-parameter deformed Virasoro algebra is presented, that is, we construct a new two-parameter quantum group related to Virasoro algebra.
Our paper is organized as follows. In section 2, we use the bosonic oscillators to construct a two-parameter deformed Virasoro algebra.
In Section 3, we present a nontrivial non-commutative and cocommutative Hopf structure of the two-parameter deformed Virasoro algebra. \vskip7pt

\cl{\bf\S2. \ The two parameters deformed Virasoro algebra $V_{p,q}$}
\setcounter{section}{2}\setcounter{theo}{0}\setcounter{equation}{0}
In this section, we compute the enveloping algebra of the two parameters deformed Virasoro algebra.\\
Using (\ref{9.100}), and it follows by induction on $n$ that
\begin{align*}
&[a,(a^{+})^{n}]=n(a^{+})^{n-1},\ \forall n\in\Z.
\end{align*}
The generators of the form
\begin{align*}
&L_{n}\equiv(a^{+})^{n+1}a
\end{align*}
realize the centerless Virasoro Lie algebra with the following bracket:
\begin{align*}
&[L_{m},L_{n}]=(n-m)L_{m+n},\ \forall m,n\in\Z.
\end{align*}
The following is an introduction to the one parameter deformation of the Virasoro algebra.\\
Fix a $q\in\C^*$ such that $q$ is not a root of unity. Instead of equation (\ref{9.100}), we assume that
\begin{align}
&[a,a^{+}]_{(1,q)}=1,\label{8.1}
\end{align}
Here we use the notation
\begin{align}
&[A,B]_{(\alpha,\beta)}=\alpha AB-\beta BA.\label{8.2}
\end{align}
From (\ref{8.1}), it follows by induction on $n$ that
\begin{align}
&[a,(a^{+})^{n}]_{(1,q^{n})}=\{n\}_{q}(a^{+})^{n-1},\label{8.3}
\end{align}
for arbitrary $n$, where the general natation
\begin{align}
&\{n\}_{q}=\frac{q^{n}-1}{q-1}\label{8.4}
\end{align}
is used. The steps of induction are as follows.
\begin{align*}
[a,(a^{+})^{2}]_{(1,q^{2})}
&=a(a^{+})^{2}-q^{2}(a^{+})^{2}a\\
&=a(a^{+})^{2}-qa^{+}aa^{+}+qa^{+}aa^{+}-q^{2}(a^{+})^{2}a\\
&=\{2\}_{q}a^{+}.
\end{align*}
Let's say that $n-1$ is true, we have
\begin{align*}
&[a,(a^{+})^{n-1}]_{(1,q^{n-1})}=\{n-1\}_{q}(a^{+})^{n-2},
\end{align*}
then
\begin{align*}
[a,(a^{+})^{n}]_{(1,q^{n})}
&=a(a^{+})^{n}-q^{n}(a^{+})^{n}a\\
&=(a(a^{+})^{n-1}-q^{n-1}(a^{+})^{n-1}a)a^{+}\\
&+q^{n-1}(a^{+})^{n-1}aa^{+}-q^{n}(a^{+})^{n}a\\
&=\{n-1\}_{q}(a^{+})^{n-1}+q^{n-1}(a^{+})^{n-1}(aa^{+}-qa^{+}a)\\
&=\{n\}_{q}(a^{+})^{n-1},
\end{align*}
so the induction holds.
\begin{prop}
The generators $L_{n}\equiv(a^{+})^{n+1}a\  (n\in\Z)$ satisfy the following relations:
\begin{align}
[L_{n},L_{m}]_{(q^{n},q^{m})}=(\{m\}_{q}-\{n\}_{q})L_{m+n},\ \forall m,n\in\Z.\label{8.5}
\end{align}
\end{prop}
\noindent{\it \textbf{Proof.}~}
Obviously, equation (5.6) holds for $m=n$ since both sides are equal to $0$. Now assume that $n \neq m$.
\begin{align*}
[L_{n},L_{m}]_{(q^{n},q^{m})}&=q^{n}L_{n}L_{m}-q^{m}L_{m}L_{n}\\
&=q^{n}(a^{+})^{n+1}a(a^{+})^{m+1}a-q^{m}(a^{+})^{m+1}a(a^{+})^{n+1}a\\
&=q^{n}(a^{+})^{n+1}\big(\{m+1\}_{q}(a^{+})^{m}+q^{m+1}(a^{+})^{m+1}a\big)a\\
&-q^{m}(a^{+})^{m+1}\big(\{n+1\}_{q}(a^{+})^{n}+q^{n+1}(a^{+})^{n+1}a\big)a\\
&=q^{n}\{m+1\}_{q}(a^{+})^{m+n+1}a+q^{m+n+1}(a^{+})^{n+m+2}aa\\
&-q^{m}\{n+1\}_{q}(a^{+})^{m+n+1}a-q^{m+n+1}(a^{+})^{n+m+2}aa\\
&=\big(q^{n}\{m+1\}_{q}-q^{m}\{n+1\}_{q}\big)(a^{+})^{m+n+1}a\\
&=\big(\{m\}_{q}-\{n\}_{q}\big)L_{m+n}.
\end{align*}
\hfill$\Box$\vskip7pt
The following is an introduction to the two parameter deformation of the Virasoro algebra. \\
  We assume that
\begin{align}
&[a,a^{+}]_{(p,q)}=paa^{+}-qa^{+}a=1.\label{8.6}
\end{align}
From (\ref{8.6}), it follows by induction on $n$ that
\begin{align}
&[a,(a^{+} )^{n}]_{(p^{n},q^{n})}=[n]_{p,q}(a^{+})^{n-1},\label{8.7}
\end{align}
for arbitrary $n$, where the quantum integer
\begin{align}
&[n]_{p,q}=\frac{p^{n}-q^{n}}{p-q}\label{8.8}
\end{align}
is used. The steps of induction are as follows.
\begin{align*}
[a,(a^{+})^{2}]_{(p^{2},q^{2})}
&=p^{2}a(a^{+})^{2}-q^{2}(a^{+})^{2}a\\
&=p\big(paa^{+}-qa^{+}a\big)a^{+}+pqa^{+}aa^{+}-q^{2}(a^{+})^{2}a\\
&=[2]_{p,q}a^{+}.
\end{align*}
Let's say that $n-1$ is true, we see that
\begin{align*}
&[a,(a^{+})^{n-1}]_{(p^{n-1},q^{n-1})}=[n-1]_{p,q}(a^{+})^{n-2},
\end{align*}
then
\begin{align*}
[a,(a^{+})^{n}]_{(p^{n},q^{n})}
&=p^{n}a(a^{+})^{n}-q^{n}(a^{+})^{n}a\\
&=p(p^{n-1}a(a^{+})^{n-1}-q^{n-1}(a^{+})^{n-1}a)a^{+}\\
&+pq^{n-1}(a^{+})^{n-1}aa^{+}-q^{n}(a^{+})^{n}a\\
&=p[n-1]_{p,q}(a^{+})^{n-1}+q^{n-1}(a^{+})^{n-1}(paa^{+}-qa^{+}a)\\
&=[n]_{p,q}(a^{+})^{n-1},
\end{align*}
so the induction holds.\\
Next we compute the enveloping algebra of the two parameters deformed Virasoro algebra.
\begin{lemm}\label{1}
The generators $L_{n}\equiv(a^{+})^{n+1}a\ (n\in\Z)$ satisfy the following relations:
\begin{align}
[L_{n},L_{m}]
&=\big(\frac{[m]_{p,q}}{p^{m}}-\frac{[n]_{p,q}}{p^{n}}\big)L_{m+n},\notag\\ &=\big(q^{n}[m]_{p,q}-q^{m}[n]_{p,q}\big)p^{-(m+n)}L_{m+n}.\label{8.9}\
\forall m,n\in\Z.
\end{align}
\end{lemm}
When $p=1$, the enveloping algebra of the two parameters deformed Virasoro algebra must be identical to the one-parameter deformation of the Virasoro algebra, we get:
\begin{align*}
[L_{n},L_{m}]
&=Aq^{n}L_{n}L_{m}-Bq^{m}L_{m}L_{n}\\
&=Aq^{n}(a^{+})^{n+1}a(a^{+})^{m+1}a-Bq^{m}(a^{+})^{m+1}a(a^{+})^{n+1}a\\
&=Aq^{n}(a^{+})^{n+1}\big(p^{-m-1}[m+1]_{p,q}(a^{+})^{m}+p^{-m-1}q^{m+1}(a^{+})^{m+1}a\big)a\\
&-Bq^{m}(a^{+})^{m+1}\big(p^{-n-1}[n+1]_{p,q}(a^{+})^{n}+p^{-n-1}q^{n+1}(a^{+})^{n+1}a\big)a,
\end{align*}
then
\begin{align*}
Aq^{n}p^{-m-1}[m+1]_{p,q}(a^{+})^{m+n+1}a-Bq^{m}p^{-n-1}[n+1]_{p,q}(a^{+})^{m+n+1}a\\
=\big(q^{n}[m]_{p,q}-q^{m}[n]_{p,q}\big)p^{-(m+n)}(a^{+})^{m+n+1}a.
\end{align*}
Because of
\begin{align*}
q^{n}[m]_{p,q}-q^{m}[n]_{p,q}=p^{-1}q^{n}[m+1]_{p,q}-p^{-1}q^{m}[n+1]_{p,q},
\end{align*}
we have
\begin{align*}
Aq^{n}p^{-m-1}[m+1]_{p,q}-Bq^{m}p^{-n-1}[n+1]_{p,q}=\big(p^{-1}q^{n}[m+1]_{p,q}-p^{-1}q^{m}[n+1]_{p,q}\big)p^{-(n+m)},
\end{align*}
then, we see that
\begin{align*}
A=p^{-n},\ B=p^{-m}.
\end{align*}
Thus, we get
\begin{align*}
[L_{n},L_{m}]=p^{-n}q^{n}L_{n}L_{m}-p^{-m}q^{m}L_{m}L_{n}.
\end{align*}
The two parameters deformed Virasoro algebra has the generating set $\{L_{n},C\mid n\in\Z\}$ and the following relations:
\begin{align}
p^{-n}q^{n}L_{n}L_{m}-p^{-m}q^{m}L_{m}L_{n}
&=\big(\frac{[m]_{p,q}}{p^{m}}-\frac{[n]_{p,q}}{p^{n}}\big)L_{m+n} \notag\\
&+\delta_{m+n,0}\frac{(q/p)^{-n}}{6\big(1+(q/p)^{n}\big)}\frac{[n-1]_{p,q}}{p^{n-1}}
\frac{[n]_{p,q}}{p^{n}}\frac{[n+1]_{p,q}}{p^{n+1}}C,\label{8.10}\\
q^{n}L_{n}C&=CL_{n}.\label{8.11}
\end{align}

\cl{\bf\S3. \ Quantum Group Structures of $V_{p,q}$}
\setcounter{section}{3}\setcounter{theo}{0}\setcounter{equation}{0}
In this section, we give a direct construction of the Hopf algebraic structures of the $V_{p,q}$.\\
The $\mathcal{U}_{p,q}$ is defined as the associative algebra generated by $\mathcal{T}$, $\mathcal{T}^{-1}$, $L_{n}$\ $(n\in\Z)$, C and relations as follows:\\
(R1) $\mathcal{T}\mathcal{T}^{-1}=1=\mathcal{T}^{-1}\mathcal{T}$;\\
(R2)
$\mathcal{T}^{m}L_{n}=p^{m(n+1)}q^{-m(n+1)}L_{n}\mathcal{T}^{m}$;\\
(R3)
$q^{m}\mathcal{T}^{m}C=p^{m}C\mathcal{T}^{m}$;\\
(R4)
$q^{n}p^{-n}L_{n}L_{m}-q^{m}p^{-m}L_{m}L_{n}=[L_{n},L_{m}]$;\\
(R5)
$q^{n}L_{n}C=p^{n}CL_{n}$;\\
Before giving the construction of the Hopf algebraic structures on $\mathcal{U}_{p,q}$, we have to check whether or not these five relations (R1)-(R5) above ensure a nontrivial associative algebra $\mathcal{U}_{p,q}$.\\
The following proposition fives a positive answer.
\begin{prop}
The associate algebra $\mathcal{U}_{p,q}$ with generators $\mathcal{T}$, $\mathcal{T}^{-1}$, $L_{n}$\ $(n\in\Z)$, C and relations (R1)-(R5) is nontrivial.
\end{prop}
\noindent{\it \textbf{Proof.}~}
Set $M:=\{L_{n},C,\mathcal{T},\mathcal{T}^{-1}\mid n\in\Z\}$. Let $T(M)$ be the tensor algebra of $M$, which is a free associative algebra generated by $M$. Then one has
\begin{align*}
T(M)=\bigoplus^{\infty}_{m=0}T(M)_{m},
\end{align*}
where $T(M)_{m}=M\otimes\ldots\otimes M$ = $span\{v_{1}\otimes\ldots\otimes v_{m}|v_{i}\in M, i = 1,\ldots,m\}.$ In particular,
\begin{align*}
T(M)_{0} = \mathbb{C}\ \ and\ \ T(M)_{1} = M.
\end{align*}
The product on $T(M)$ is naturally defined by \begin{align*}
(v_{1}\otimes\ldots\otimes v_{m})(w_{1}\otimes\ldots\otimes w_{m})=v_{1}\otimes\ldots\otimes v_{m}\otimes w_{1}\otimes\ldots\otimes w_{m}.
\end{align*}
Let I be the two-sided ideal of $T(M)$ generated by
\begin{align}
&\mathcal{T}\otimes\mathcal{T}^{-1}-\mathcal{T}^{-1}\otimes\mathcal{T};\label{9.1}\\
&T^{m}\otimes L_{n}-p^{m+n}q^{-m(n+1)}L_{n}\otimes T^{m};\label{9.2}\\
&q^{n}p^{-(m+n)}L_{n}\otimes L_{m}-q^{m}p^{-(m+n)}L_{m}\otimes L_{n};\label{9.3}\\
&q^{m}T^{m}\otimes C-C\otimes T^{m}\label{9.4}
\end{align}
for all $m,n\in\Z$ and where $T^{-n}=(T^{-1})^{n}$. Set
\begin{align*}
S(M):=T(M)/I.
\end{align*}
It is obvious that S(M) is also a $\mathbb{Z}$-graded algebra with a basis
\begin{align}
\widetilde{B}=\{T^{d}(\mathcal{T}^{-1})^{d^{'}}L^{k_{i_{1}}}_{i_{1}}\ldots L^{k_{i_{m}}}_{i_{m}}C^{d_{1}}\},\label{9.5}
\end{align}
where $k_{i_{p}},d,d^{'},d_{1}\in\mathbb{N}$; $i_{p}\in\mathbb{Z}_{+}$; $i_{1}<\ldots<i_{m}$. Let $\widetilde{J}$ be another two-sided ideal of $T(M)$ generated by the elements form
\begin{align}
q^{n}p^{-(m+n)}L_{n}\otimes L_{m}-q^{m}p^{-(m+n)}L_{m}\otimes L_{n}-[L_{n},L_{m}]\label{9.6}
\end{align}
together with those in (\ref{9.1}), (\ref{9.2}) and (\ref{9.4}). Then set
\begin{align*}
\widetilde{\mathcal{U}_{p,q}}:=T(M)/\widetilde{J}.
\end{align*}
Our aim is to show that $\widetilde{B}$ defined in (\ref{9.5}) is also a basis of $\widetilde{\mathcal{U}_{p,q}}$. Let
\begin{align*}
\widetilde{B}^{'}=\{v_{i_{1}}\otimes\ldots\otimes v_{i_{m}}|v_{i}\in M,1\leq i_{1}\leq i_{2}\leq\ldots\leq i_{m},m\geq0\}
\end{align*}
be a subset of $T(M)$ and let $U^{'}$  be the subspace of $T(M)$ spanned by $\widetilde{B}^{'}$. We claim that
\begin{align}
T(M)=U^{'}\oplus \widetilde{J}.\label{9.7}
\end{align}
For any $v\in T(M)$, we can write $v=v^{(m)}+v^{(m-1)}+\ldots+v^{(0)}$, where $v^{(m)}\neq0$ for some $m\geq0$ and where $v^{(i)}\in T(M)_{i}$ with $i=0,1,\ldots,m$. We call $m$ the degree of $v$. From (\ref{9.1}), (\ref{9.2}), (\ref{9.4}) and (\ref{9.6}), it follows
\begin{align*}
v_{i_{1}}\otimes\ldots\otimes\big(v_{i_{k}}\otimes v_{i_{k+1}}-v_{i_{k+1}}\otimes v_{i_{k}}-[v_{i_{k}},v_{i_{k+1}}]\big)\otimes\ldots v_{i_{m}}\in \widetilde{J},
\end{align*}
namely, the difference between $v_{i_{1}}\otimes\ldots\otimes v_{i_{k}}\otimes v_{i_{k+1}}\otimes\ldots\otimes v_{i_{m}}$
and $av_{i_{1}}\otimes\ldots\otimes v_{i_{k+1}}\otimes v_{i_{k}}\otimes\ldots\otimes v_{i_{m}}$ (for some $a\in \mathbb{C}^{*}$) is an element in $\widetilde{J}$ and an element with degree than $m$. So by induction on the degree of $v$ one can obtain that $T(M)=U^{'}+\widetilde{J}$.\\
\ \ It remains to show that equation (\ref{9.7}) is a direct sum, which is equivalent to the linear independence of $\widetilde{B}$ in $\widetilde{\mathcal{U}_{p,q}}$. Suppose that a nonzero linear combination $v$ of the elements in $\widetilde{B}^{'}$ is in $\widetilde{J}$. It follows from (\ref{9.1}), (\ref{9.2}), (\ref{9.4}) and (\ref{9.6}) that homogeneous component $v^{(m)}$ of $v$ with highest degree must lie in ker$\pi$ (by comparing (\ref{9.3}) and (\ref{9.6})), where $\pi$:
$T(M)\rightarrow S(M)$ is the natural $\mathbb{Z}$-graded algebraic homomorphism, namely,
\begin{align*}
\pi(v_{i_{1}}\otimes\ldots\otimes v_{i_{m}})=v_{i_{1}}v_{i_{2}}\ldots v_{i_{m}}.
\end{align*}
However, $v^{(m)}$ is a nonzero linear combination of the elements in $\widetilde{B}^{'}$, it is impossible to appear in ker$\pi$. This contradiction implies $\widetilde{B}$ is a basis of $\widetilde{\mathcal{U}_{p,q}}$. Since it is clear that $\mathcal{U}_{p,q}\cong \widetilde{\mathcal{U}_{p,q}}/J$, where $J$ is the two-sided ideal of $\widetilde{\mathcal{U}_{p,q}}$ generated by $\mathcal{T}\mathcal{T}^{-1}-1$, we obtain taht nontrivial associative algebra with basis
\begin{align}
\widetilde{B}^{'}=\{T^{d}L^{k_{i_{1}}}_{i_{1}}\ldots L^{k_{i_{m}}}_{i_{m}}C^{d_{1}}\}\label{9.8}
\end{align}
where $d\in\mathbb{Z}$;
$d_{1}\in\mathbb{N}$;
$i_{p},j_{q}\in\mathbb{Z}_{+}(i,p=1,2,\ldots,m;j,q=1,2,\ldots,n)$; $i_{1}<\ldots<i_{m}$.
\hfill$\Box$\vskip7pt
With the above proposition in hand, we can safely proceed with the construction of the Hopf algebraic structures on $\widetilde{\mathcal{U}_{p,q}}$ now. This will be done by severall lemmas below.
\begin{lemm}
There is a unique algebraic homomorphism $\Delta$: $\mathcal{U}_{p,q}\rightarrow \mathcal{U}_{p,q}\times \mathcal{U}_{p,q}$ with
\begin{align}
&\Delta(\mathcal{T})=\mathcal{T}\otimes \mathcal{T},
\Delta(\mathcal{T}^{-1})=\mathcal{T}^{-1}\otimes \mathcal{T}^{-1},\label{9.9}\\
&\Delta(L_{n})=L_{n}\otimes \mathcal{T}^{n}+\mathcal{T}^{n}\otimes L_{n},\label{9.10}\\
&\Delta(C)=C\otimes 1+1\otimes T,\label{9.11}\\
&\epsilon(L_{n})=\epsilon(C)=0,\label{9.12}\\
&S(L_{n})=-\mathcal{T}^{-n}L_{n}\mathcal{T}^{-n},\label{9.13}\\
&S(C)=-C.\label{9.14}
\end{align}
\end{lemm}
\noindent{\it \textbf{Proof.}~}
It is clear that $\Delta(\mathcal{T}^{m})=\mathcal{T}^{m}\otimes \mathcal{T}^{m}$ for arbitrary $m\in\mathbb{Z}$, $\epsilon(\mathcal{T}^{m})=1$, $S(\mathcal{T}^{m})=\mathcal{T}^{-m}$. We see that $\Delta(\mathcal{T})$, $\Delta(\mathcal{T}^{-1})$, $\Delta(L_{n})$ satisfy the relations (R1)-(R5). This is trivial for (R1). For (R2) and (R3) it follows directly from (\ref{9.9})-(\ref{9.14}).Now look at (R4), we see that
\begin{align*}
\Delta(L_{n})\Delta(L_{m})&=(L_{n}\otimes \mathcal{T}^{n}+\mathcal{T}^{n}L_{n})(L_{m}\otimes \mathcal{T}^{m}+\mathcal{T}^{m}L_{m})\\
&=L_{n}L_{m}\otimes \mathcal{T}^{m+n}+L_{n}\mathcal{T}^{m}\otimes \mathcal{T}^{n}L_{m}+\mathcal{T}^{n}L_{m}\otimes L_{n}\otimes \mathcal{T}^{m}+\mathcal{T}^{m+n}\otimes L_{n}L_{m}\\
&=L_{n}L_{m}\otimes \mathcal{T}^{m+n}+p^{n(m+1)}q^{-n(m+1)}L_{n}\mathcal{T}^{m}\otimes L_{m}\mathcal{T}^{n}\\
&+p^{n(m+1)}q^{-n(m+1)}L_{m}\mathcal{T}^{n}\otimes L_{n}\mathcal{T}^{m}
+\mathcal{T}^{m+n}\otimes L_{n}L_{m}.
\end{align*}
Similarly, we get
\begin{align*}
\Delta(L_{m})\Delta(L_{n})
&=L_{m}L_{n}\otimes \mathcal{T}^{m+n}+p^{m(n+1)}q^{-m(n+1)}L_{m}\mathcal{T}^{n}\otimes L_{n}\mathcal{T}^{m}\\
&+p^{m(n+1)}q^{-m(n+1)}L_{n}\mathcal{T}^{m}\otimes L_{m}\mathcal{T}^{n}
+\mathcal{T}^{m+n}\otimes L_{m}L_{n}.
\end{align*}
Then it follows
\begin{align*}
&q^{n}p^{-n}\Delta(L_{n})\Delta(L_{m})-
q^{m}p^{-m}\Delta(L_{m})\Delta(L_{n})\\
&=(q^{n}p^{-n}L_{n}L_{m}\otimes \mathcal{T}^{m+n}-q^{m}p^{-m}L_{m}L_{n}\otimes \mathcal{T}^{m+n})+(q^{n}p^{-n}\mathcal{T}^{m+n}\otimes L_{n}L_{m}-q^{m}p^{-m}\mathcal{T}^{m+n}\otimes L_{m}L_{n})\\
&=\big(\frac{[m]_{p,q}}{p^{m}}-\frac{[n]_{p,q}}{p^{n}}\big)
(L_{m+n}\otimes \mathcal{T}^{m+n}+\delta_{m+n,0}\frac{(q/p)^{-n}}{6\big(1+(q/p)^{n}\big)}\frac{[n-1]_{p,q}}{p^{n-1}}
\frac{[n]_{p,q}}{p^{n}}\frac{[n+1]_{p,q}}{p^{n+1}}C\otimes \mathcal{T}^{m+n})\\
&+\big(\frac{[m]_{p,q}}{p^{m}}-\frac{[n]_{p,q}}{p^{n}}\big)
(\mathcal{T}^{m+n}\otimes L_{m+n}+\mathcal{T}^{m+n}\otimes \delta_{m+n,0}\frac{(q/p)^{-n}}{6\big(1+(q/p)^{n}\big)}\frac{[n-1]_{p,q}}{p^{n-1}}
\frac{[n]_{p,q}}{p^{n}}\frac{[n+1]_{p,q}}{p^{n+1}}C)\\
&=\big(\frac{[m]_{p,q}}{p^{m}}-\frac{[n]_{p,q}}{p^{n}}\big)
\Delta(L_{m+n})+\delta_{m+n,0}\frac{(q/p)^{-n}}{6\big(1+(q/p)^{n}\big)}\frac{[n-1]_{p,q}}{p^{n-1}}
\frac{[n]_{p,q}}{p^{n}}\frac{[n+1]_{p,q}}{p^{n+1}}\Delta(C)
\end{align*}
In addition, we need to show that $(S(\mathcal{T}),S(\mathcal{T}^{-1}),S(L_{m}))$ satisfies the relations (R1)-(R5) in $\mathcal{U}_{p,q}^{opp}$. Let us denote the multiplication in $\mathcal{U}_{p,q}^{opp}$ by a $"\cdot"$ in order to distinguish it from that in $\mathcal{U}_{p,q}$. It is easy to get $S(\mathcal{T}^{m})=\mathcal{T}^{-m}\ (m\in\mathbb{Z})$, For (R2), we get:\\
\begin{align*}
S(\mathcal{T}^{m})\cdot S(L_{m})
&=S(L_{m})S(\mathcal{T}^{m})\\
&=-\mathcal{T}^{-m}L_{m}\mathcal{T}^{-m}\mathcal{T}^{-m}\\
&=-\mathcal{T}^{-m}p^{m(m+1)}q^{-m(m+1)}\mathcal{T}^{-m}L_{m}\mathcal{T}^{-m}\\
&=p^{m(m+1)}q^{-m(m+1)}S(\mathcal{T}^{m})S(L_{m})\\
&=p^{m(m+1)}q^{-m(m+1)}S(L_{m})\cdot S(\mathcal{T}^{m}).
\end{align*}
For (R4):
\begin{align*}
q^{n}p^{-n}S(L_{n})\cdot S(L_{m})
&=q^{n}p^{-n}\mathcal{T}^{-m}L_{m}\mathcal{T}^{-m}\mathcal{T}^{-n}L_{n}\mathcal{T}^{-n}\\
&=q^{n}p^{-n}p^{n(m+1)}q^{-n(m+1)}\mathcal{T}^{-m-n}L_{m}p^{-m(n+1)}q^{m(n+1)}L_{n}\mathcal{T}^{-m}L_{n}\mathcal{T}^{-n}\\
&=q^{m}p^{-m}\mathcal{T}^{-m-n}L_{m}L_{n}\mathcal{T}^{-m-n},
\end{align*}
thus
\begin{align*}
q^{m}p^{-m}S(L_{m})\cdot S(L_{n})=q^{n}p^{-n}\mathcal{T}^{-m-n}L_{n}L_{m}\mathcal{T}^{-m-n}.
\end{align*}
Then we see that
\begin{align*}
&q^{n}p^{-n}S(L_{n})\cdot S(L_{m})-q^{m}p^{-m}S(L_{m})\cdot S(L_{n})\\
&=\mathcal{T}^{-m-n}\Big(-\big(\frac{[m]_{p,q}}{p^{m}}-\frac{[n]_{p,q}}{p^{n}}
\big) -\delta_{m+n,0}\frac{(q/p)^{-n}}{6\big(1+(q/p)^{n}\big)}\frac{[n-1]_{p,q}}{p^{n-1}}
\frac{[n]_{p,q}}{p^{n}}\frac{[n+1]_{p,q}}{p^{n+1}}C
  \Big)\mathcal{T}^{-m-n}\\
&=\big(\frac{[m]_{p,q}}{p^{m}}-\frac{[n]_{p,q}}{p^{n}}
\big) S(L_{m+n})+\delta_{m+n,0}\frac{(q/p)^{-n}}{6\big(1+(q/p)^{n}\big)}\frac{[n-1]_{p,q}}{p^{n-1}}
\frac{[n]_{p,q}}{p^{n}}\frac{[n+1]_{p,q}}{p^{n+1}}S(C).
\end{align*}
One can similarly check that (R3) and (R5) are also preserved by $S$. So there is indeed a homomorphism $S:\mathcal{U}_{p,q}\rightarrow \mathcal{U}_{p,q}^{opp}$ or an antihomomorphism $S:\mathcal{U}_{p,q}\rightarrow \mathcal{U}_{p,q}$ satisfying (\ref{9.13}) and (\ref{9.14}). Now $S^{2}$ is an ordinary homomorphism from $\mathcal{U}_{p,q}$ to $\mathcal{U}_{p,q}$.
One can check easily on the  generators that $S^{2}=id$, which implies that $S$ is bijective. It is straightforward to see that $\epsilon(\mathcal{T})=\epsilon(\mathcal{T}^{-1})=1, \epsilon(L_{n})=\epsilon(C)=0$ satisfy the relations (R1)-(R5). So we have the algebraic homomorphism $\epsilon$.
\hfill$\Box$\vskip7pt
\begin{theo}
$(\mathcal{U}_{p,q},\Delta,\epsilon,S)$ defined by (R1)-(R5) and (\ref{9.9})-(\ref{9.14}) is a Hopf algebra.
\end{theo}
\textbf{Corollary 5.7}
As vector spaces, one has
\begin{align*}
\mathcal{U}_{p,q}\cong \mathbb{C}[\mathcal{T},\mathcal{T}^{-1}]\otimes_{\mathbb{C}}U_{p,q},
\end{align*}
where $U_{p,q}=U(V_{p,q})$ is the enveloping algebra of $V_{p,q}$ generated by $L_{n}(n\in \mathbb{Z})$ with relations (\ref{8.5}).\\

\end{document}